\documentclass[journal]{IEEEtran}

\usepackage{amsmath}
\usepackage{amssymb,amsfonts}
\usepackage{graphicx}
    \graphicspath{{./pix/}}
\usepackage{subfigure} 
\usepackage[all]{xy}
\usepackage{epstopdf}

\begin{document}
\title{
Data Assimilation by Conditioning on Future Observations
}

\author{
Wonjung Lee and Chris Farmer
\thanks{
Wonjung Lee and Chris Farmer are with 
the Oxford Centre for Collaborative Applied Mathematics
(OCCAM)
in
the Mathematical Institute,
University of Oxford, Oxford, U.K. 
(email:leew@maths.ox.ac.uk; farmer@maths.ox.ac.uk).}
}

\maketitle

\begin{abstract}
Conventional recursive filtering approaches, designed for quantifying the state of an evolving uncertain
dynamical system
with intermittent observations,
use a sequence of (i) an uncertainty propagation step followed by (ii) a step where
the associated data is
assimilated using Bayes' rule.
In this paper we switch the order of the steps to: (i) one step ahead data assimilation
followed by (ii) uncertainty propagation. This route leads to a class of filtering algorithms 
named \emph{smoothing filters}.
For a system driven by random noise, our proposed methods require the 
probability distribution
of the driving noise 
after the assimilation
to be biased 
by a nonzero mean.
The system noise, conditioned on future observations,
in turn pushes forward the filtering solution in time closer to the true state
and indeed helps to find a more accurate 
approximate solution for the state estimation problem.

\end{abstract}


\IEEEpeerreviewmaketitle

\section{Introduction}
\noindent
\IEEEPARstart{T}{here} are many problems
in science and engineering in which the state of a system has to be identified 
from a set of noisy observations.
The solution has concrete applications 
in fields such as statistical signal processing, 
sonar ranging, target tracking, satellite navigation,
and prediction of weather and climate in atmosphere-ocean dynamics
\cite{
jazwinski1970stochastic, 
anderson1979optimal, 
doucet2001sequential,
evensen2009data}.

Many of these problems involve
(i) a forward model for the state evolution of a dynamical system
and (ii) observational data associated with the system state.
In filtering, one combines these two process equations
to form an effective solution of the state estimation problem.
For practical reasons, online estimation via a recursive method is desired.
The conventional way to achieve this real-time filtering algorithm
is to alternate application of (i) the uncertainty quantification (UQ)
or time update and (ii) the data assimilation (DA)
or measurement update, in a sequential fashion.
The UQ corresponds to solving the Fokker-Planck
equation for continuous-time dynamical systems, 
and the Chapman-Kolmogorov equation for discrete-time dynamical systems.
The algorithm, called a Bayesian filter, achieves DA using Bayes' rule.

For linear dynamics with a linear observation process, 
the filtering solution is Gaussian and the Kalman
filter provides the answer
\cite{Kalman60}.
When nonlinearity is present in either the forward model or the observation process, 
one in general has to develop an approximate solution due to the lack of an analytical solution.
One such approximation is the extended Kalman filter
based on successive linearisation of both the forward model
and the 
observation process
\cite{gelb1974applied}.
Other examples of Gaussian approximation filters are the unscented Kalman filters
\cite{julier2004unscented},
cubature Kalman filters
\cite{arasaratnam2009cubature}
and Gaussian particle filters
\cite{kotecha2003gaussian}.
In place of linearisation, these filters use 
discrete measures
for a Gaussian approximation
of the conditioned probability distribution.
There are other filters, called Gaussian sum filters,
where the filtering solution is
represented by multiple weighted Gaussian kernels rather than a single one
(see for example
\cite{chen2000mixture,
stordal2011bridging}).
The ensemble Kalman filter 
\cite{evensen2009data}
and
the bootstrap filter 
\cite{gordon1993novel}
are sequential Monte Carlo methods where discrete measures are used to 
approximate the probability distribution.
These two filters
can be viewed as specific cases of 
one kind of Gaussian sum filter
\cite{stordal2011bridging}.

\section{Sequential Data Assimilation}
\label{sec:bayesianfilter}
\noindent
Let the discrete-time evolution of an $\mathbb{R}^{d}$-valued vector, $\mathbf{x}$, be  
governed by the
\begin{equation}
  \label{eq:fmodel}
  \textbf{Forward model} \quad 
  \mathbf{x}_{n+1}=\Phi^n( \mathbf{x}_{n} , {\xi}_n), 
  \quad \xi_n \sim \mathcal{N}\left( \mathbf{0}, \Gamma_n \right)
\end{equation}
where 
$n \in \mathbb{N} \cup \{0\}$ labels the time step,
$\xi_n \in \mathbb{R}^{D}$ is an 
independent and identically distributed (i.i.d.) Gaussian noise
and $\mathbf{0}$ denotes a zero vector (or later, a zero matrix).
Data $\mathbf{y}_{n} \in \mathbb{R}^{d'}$, 
associated with $\mathbf{x}_n$, is modelled by the
\begin{equation}
  \label{eq:observation}
  \textbf{Observation} \quad \mathbf{y}_{n} = 
  \phi^n(\mathbf{x}_{n}) +\eta_{n}, \quad \eta_{n} \sim \mathcal{N}(\mathbf{0},R_{n})
\end{equation}
for a measurement function $\phi^n$
and i.i.d. Gaussian $\eta_n$.
Here 
$\mathbf{x}_n$, $\xi_n$ and $\eta_n$ are assumed statistically independent.
Let
$
\mathcal{X}_n \equiv 
\left[
\begin{array}{c}
  \mathbf{x}_n  \\ 
  \xi_n
\end{array}\right]
$
be the $(d+D)$-dimensional augmented system.
Let $\mathbf{x}_{n|n'} \equiv \mathbf{x}_n|Y_{n'}$,
$\xi_{n|n'} \equiv \xi_n|Y_{n'}$ and
$\mathcal{X}_{n|n'} \equiv \mathcal{X}_n|Y_{n'}$
be random vectors
conditioned on the collection of observations 
$Y_{n'} \equiv \lbrace \mathbf{y}_1,\cdots,\mathbf{y}_{n'} \rbrace$.
It is called \emph{smoothing} 
to find the probability distribution of $\mathbf{x}_{n|n'}$ or $\mathcal{X}_{n|n'}$
when $n<n'$, \emph{filtering} when $n=n'$, 
and \emph{prediction} when $n>n'$.

\subsection{Conventional Filters and Smoothing Filters}
Given Eqs.~(\ref{eq:fmodel}), (\ref{eq:observation})
and the probability distribution of the initial condition $\mathbf{x}_0$,
the sequential filtering problem requires
finding the probability distribution of
$\mathbf{x}_{n|n}$ for $n \geq 1$.
The conventional approach
to such problems
is to alternate the time update
$\mathbf{x}_{n|n} \to \mathbf{x}_{n+1|n}$
or $\mathcal{X}_{n|n} \to \mathbf{x}_{n+1|n}$
for prediction
and the measurement update
$\mathbf{x}_{n+1|n} \Rightarrow \mathbf{x}_{n+1|n+1}$
for filtering,
in a sequential fashion.
We here use the notation $\to$ to increase the first index by one 
or UQ and the notation $\Rightarrow$ to increase the second index 
by one or DA.
We call a method following the conventional approach, a 
\emph{conventional filter}.
While most of the 
prevailing filters fall into this category,
it is of course possible to solve the sequential filtering problem by 
following the other route, i.e.,
successive application of a
measurement update
$\mathcal{X}_{n|n} \Rightarrow \mathcal{X}_{n|n+1}$
for smoothing
and a time update
$\mathcal{X}_{n|n+1} \to \mathbf{x}_{n+1|n+1}$
for filtering
\cite{desbouvries2011direct}.
Such an algorithm achieves the data assimilation via smoothing
and so we call the algorithm a \emph{smoothing filter}.
The following two  sequential methods:
\begin{enumerate}
  \item Conventional filter :
	$\mathcal{X}_{n|n} \to \mathbf{x}_{n+1|n} \Rightarrow \mathbf{x}_{n+1|n+1}$
  \item Smoothing filter :
	$\mathcal{X}_{n|n} \Rightarrow \mathcal{X}_{n|n+1} \to \mathbf{x}_{n+1|n+1}$
\end{enumerate}
are illustrated in 
Fig.~\ref{fig:filter}.

\begin{figure}[t!]
\begin{center}
\begin{displaymath}
    \xymatrix{
	  \mathcal{X}_{n|n} = \left[ \mathbf{x}_{n|n}, \xi_{n|n} \right]
	  \ar[r]^{} \ar@2[d]_{} 
	  & \mathbf{x}_{n+1|n} \ar@2[d]^{} \\
	  \mathcal{X}_{n|n+1} 
	  \ar[r]_{} 
	  & \mathbf{x}_{n+1|n+1} }
\end{displaymath}
\end{center}
	\caption[]{
	  Conventional filter and smoothing filter.
	}
\label{fig:filter}
\end{figure}
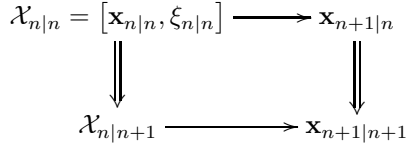

\begin{figure}[t!]
\begin{center}
\begin{displaymath}
    \xymatrix{
	  \mathcal{X}_{n|n}  \ar[r]^{\mathbf{y}_{n+1}} \ar@2[d]_{} 
	  & \widehat{\mathbf{x}}_{n+1} \ar@{|->}[r]^{}
	  & \mathbf{x}_{n+1|n} \ar@2[d]^{} \\
	  \mathcal{X}_{n|n+1}
	  \ar[rr]^{\xi_{n}\vert \mathbf{y}_{n+1}} 
	  &
	  & \mathbf{x}_{n+1|n+1} }
\end{displaymath}
\end{center}
	\caption[]{
	  Importance sampling characteristic of smoothing filter.
	}
\label{fig:filter2}
\end{figure}
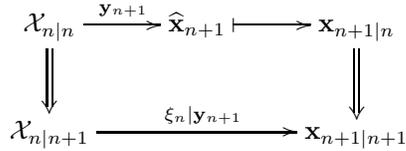

\subsection{
 An automatic form of importance sampling
}
\label{subsec:Importance}
A naive conventional filter can unfortunately be quite in error
when the distance between
$\mathbf{y}_{n+1}$ 
and
the approximation of $\mathbf{x}_{n+1|n}$
is large in some sense,
and therefore the application of such a
DA algorithm fails to produce an accurate approximation of $\mathbf{x}_{n+1|n+1}$.
One approach for reducing this distance, 
with the goal of improving the approximation of $\mathbf{x}_{n+1|n+1}$,
introduces a new forward model
\begin{equation}
  \label{eq:fmodel3}
  \widehat{\mathbf{x}}_{n+1} = \widehat{\Phi}^n(\widehat{\mathbf{x}}_{n}, \xi_n ; \mathbf{y}_{n+1})
\end{equation}
depending on the value of $\mathbf{y}_{n+1}$.
An appropriate choice of the operator 
$\widehat{\Phi}^n$
in Eq.~(\ref{eq:fmodel3})
with $\widehat{\mathbf{x}}_n = \mathbf{x}_{n|n}$,
if combined 
with a subsequent mapping 
$\widehat{\mathbf{x}}_{n+1} \mapsto \mathbf{x}_{n+1|n}$,
would lead to
an approximation 
$\mathbf{x}_{n+1|n}$
closer to $\mathbf{y}_{n+1}$. 
In this way,
the two step approximation
through the intermediate variable $\widehat{\mathbf{x}}_{n+1}$
can produce an accurate filtering solution.
The method of 
sequential importance sampling 
is one kind of Monte Carlo bootstrap filters 
and
is developed following this idea 
\cite{doucet2000sequential,
van2010nonlinear}.

It is worth remarking that a smoothing filter
can be a competitive approach due to
its inherent importance sampling characteristic.
More precisely, in contrast to the case of conventional filters for which
$\xi_{n|n} = \xi_n$ and
$\mathbb{E}(\xi_{n|n}) = \mathbf{0}$,
the conditioned variable $\xi_{n|n+1}$
in a smoothing filter 
is biased in the sense
$\mathbb{E}(\xi_{n|n+1}) \neq \mathbf{0}$,
where 
$\mathbb{E}(  \cdot )$ denotes the statistical average.
Note that, though the observation $\mathbf{y}_{n+1}$ is not directly involved
in the UQ of smoothing filters,
the data assimilated 
driving noise
$\xi_{n|n+1}$ 
can give a nontrivial effect to 
the UQ from
$\mathbf{x}_{n+1|n+1}=\Phi^n( \mathbf{x}_{n|n+1} , {\xi}_{n|n+1})$
and make
the approximation of $\mathbf{x}_{n+1|n+1}$ 
closer to $\mathbf{y}_{n+1}$,
compared with the one from the conventional filter
(see Fig.~\ref{fig:filter2}).
In other words,
without introducing an additional forward model,
smoothing filters achieve a similar effect with
Monte Carlo importance sampling.
It becomes therefore our aim to develop 
smoothing filters
and investigate their possible outperformances in solving the state estimation problem.

\subsection{Gaussian approximation filtering}
One way of building a smoothing filter is to mimic
an existing conventional filter.
In this paper, 
we develop a smoothing filter
by letting its UQ and DA methods
be basically of the same kind
as those of a conventional filter,
but with the ordering of UQ and DA reversed.

Let us note that there are many conventional filters
that adopt a sum of  Gaussian kernels
(or Dirac masses)
to approximate the conditioned probability distribution,
and although these may ultimately be the approximations of choice, in this paper concern is confined to the problem of approximation using a single Gaussian density.
This is because the aim is to develop simple but efficient filtering methods, similar to traditional methods, but with enhanced accuracy. However, we hope that any success in this aim will help direct future efforts toward the wider aim of  developing a rigorous and convergent method, perhaps using a Gaussian sum approximation.

The rest of the paper is organised as follows.
We develop a number of conventional Gaussian filters
based on some traditional filters
in Section~\ref{sec:gaussianfiltering}
and 
formulate corresponding Gaussian smoothing filters
in Section~\ref{sec:variationalsmoothing}.
With the help of the
test problems
gathered in Section~\ref{sec:comparison},
numerical simulations are used to examine the accuracy of smoothing filters 
in Section~\ref{eq:numericalsimulations}.
We conclude our results in Section~\ref{sec:conclusion}.

\section{Conventional Gaussian Filtering}
\label{sec:gaussianfiltering}
By combining the UQ methods presented in subsection~\ref{subsec:TU}
and the DA methods presented in subsection~\ref{subsec:MU},
we develop a number of conventional Gaussian approximation filters
in subsection~\ref{subsec:CGF}.
Because the mean and covariance completely determine the Gaussian distribution,
the algorithm defines the mapping of the first two moments.

Let the mean and covariance of $\mathbf{x}_{n|n'}$ be denoted by
$\bar{\mathbf{x}}_{n|n'}$ and $\mathbf{C}_{n|n'}$.
Let the mean and covariance of $\mathcal{X}_{n|n'}$ be denoted by $\bar{\mathcal{X}}_{n|n'}$ and
$\mathcal{C}_{n|n'}$.

\subsection{Time Update ($\mathcal{X}_{n|n} \to \mathbf{x}_{n+1|n}$)}
\label{subsec:TU}
From now on
we use the notation 
$\mathbf{x}_{n+1}={\Phi}^n(\mathcal{X}_n)$
in place of Eq.~(\ref{eq:fmodel}). 
When $\mathcal{X}_{n|n'}$ is Gaussian, 
the following two 
Gaussian approximations for $\mathbf{x}_{n+1|n'}=\Phi^n(\mathcal{X}_{n|n'})$
can be used to quantify the uncertainty propagation.
The first method makes use of a linearisation of $\Phi^n(\mathcal{X}_{n|n'})$
to obtain an approximate Gaussian random variable.
The second one
approximates
$\Phi^n(\mathcal{X}_{n|n'})$ itself,
which is not Gaussian unless $\Phi^n$ is a linear function,
by a Gaussian whose mean and covariance are given by
those of
$\Phi^n(\mathcal{X}_{n|n'})$.

\subsubsection{Linear Gaussian approximation}
\noindent
The first order Taylor approximation
\begin{equation*}
  \Phi^n({\mathcal{X}}_{n|n'}) \simeq \Phi^n(\bar{\mathcal{X}}_{n|n'})+
  \nabla \Phi^n \vert_{\bar{\mathcal{X}}_{n|n'}}
  (\mathcal{X}_{n|n'}-\bar{\mathcal{X}}_{n|n'})
\end{equation*}
leads to
\begin{equation}
  \begin{split}
	\label{eq:momentmapping}
	\bar{\mathbf{x}}_{n+1|n'} & \simeq \Phi^n(\bar{\mathcal{X}}_{n|n'}), \\
	\mathbf{C}_{n+1|n'} & \simeq \nabla \Phi^n \vert_{\bar{\mathcal{X}}_{n|n'}} \mathcal{C}_{n|n'}
	\; (\nabla \Phi^n \vert_{\bar{\mathcal{X}}_{n|n'}})^T
  \end{split}
\end{equation}
where the superscript $T$ denotes the matrix transpose.
An application
of Eq.~(\ref{eq:momentmapping})
with $n'=n$,
$
\bar{\mathcal{X}}_{n|n}=
\left[
\begin{array}{c}
  \bar{\mathbf{x}}_{n|n}  \\ 
  \mathbf{0}
\end{array}\right]
$
and
$
\mathcal{C}_{n|n}=\left[
\begin{array}{cc}
  \mathbf{C}_{n|n} & \mathbf{0} \\ 
  \mathbf{0} & \Gamma_n
\end{array}\right]
$ 
produces one prediction algorithm.

\subsubsection{Point-based Gaussian approximation}
Let $\delta_x$ denote a Dirac mass centered at $x$.
Let 
$\sum_{j} \Lambda_j \delta_{\mathcal{X}^j_{n|n'}}$
be a discrete measure approximating the law of $\mathcal{X}_{n|n'}$.
Then the discrete measure 
$\sum_j \Lambda_j \delta_{\Phi^n({\mathcal{X}^j_{n|n'}})}$
approximates the law of 
$\Phi^n(\mathcal{X}_{n|n'})$.
By mapping the first two moments according to the equations,
\begin{equation}
  \begin{split}
	\label{eq:momentmapping2}
&	\bar{\mathbf{x}}_{n+1|n'}  \simeq \sum_j \Lambda_j \Phi^n({\mathcal{X}^j_{n|n'}}), \\
&	\mathcal{C}_{n+1|n'}  \simeq \\
&  \quad \sum_j \Lambda_j 
	\left( \Phi^n({\mathcal{X}^j_{n|n'}}) - \bar{\mathbf{x}}_{n+1|n'} \right)
	\left( \Phi^n({\mathcal{X}^j_{n|n'}}) - \bar{\mathbf{x}}_{n+1|n'} \right)^T
  \end{split}
\end{equation}
with $n'=n$ 
one derives the point-based Gaussian approximation algorithm.
In developing the approximate measures as the point-based approximation of $\mathcal{X}_{n|n'}$, one
can use (i) cubature measure supported on deterministically placed points 
or (ii) empirical measure, i.e., an equally weighted discrete measure supported on a set of random points.
Recall that a weighted discrete measure
is called a cubature measure of degree $r$
with respect to the given probability distribution 
provided the moments 
of these two measures
agree with one another up to total degree $r$.

We mention that
Eqs.~(\ref{eq:momentmapping}) and ~(\ref{eq:momentmapping2})
are very similar to the prediction methods used in the
extended Kalman filter and 
cubature Kalman filter (cubature measure)
or
Gaussian particle filter (empirical measure),
respectively.
The algorithms in these traditional filters
compute
$\mathbf{x}_{n|n} \to \mathbf{x}_{n+1|n}$
for the forward model
\begin{equation*}
  \mathbf{x}_{n+1}=f^n( \mathbf{x}_{n} )+  {\xi}_n
\end{equation*}
which is a specific case of
Eq.~(\ref{eq:fmodel}).
We derive
the formula 
$\mathcal{X}_{n|n} \to \mathbf{x}_{n+1|n}$
because the algorithm
can be used 
(i)
for a general forward model including the case of multiplicative noise,
and
(ii)
for the smoothing filters as will be developed in the next section.

\subsection{Measurement Update ($\mathbf{x}_{n+1|n} \to \mathbf{x}_{n+1|n+1}$)}
\label{subsec:MU}
Bayes' rule,
\begin{equation}
  \label{eq:bayesrule}
  \mathbb{P}(X|Y)  = \frac{\mathbb{P}(X,Y) }{\mathbb{P}(Y) } 
\end{equation}
for random vectors $X$ and $Y$
can be employed to 
find the conditioned probability distribution.
Eq.~(\ref{eq:bayesrule}) implies that
if $X$ and $Y$ are jointly Gaussian, i.e.,
$Z = 
\left[
\begin{array}{c}
  X  \\ 
  Y 
\end{array}\right]
$
is Gaussian with mean 
$
\left[
\begin{array}{c}
  \bar{x}  \\ 
  \bar{y} 
\end{array}\right]
$
and covariance
$
\left[
\begin{array}{cc}
  \Sigma_{xx} & \Sigma_{xy} \\ 
  \Sigma_{yx} & \Sigma_{yy}
\end{array}\right]
$,
then 
the conditioned variable 
$X\vert Y$ 
with
$Y=y$ 
is Gaussian with
mean and covariance
given by
\begin{equation}
  \label{eq:conditionedmoments}
  \begin{split}
	\bar{x}'  & = \bar{x}+\Sigma_{xy}\Sigma_{yy}^{-1}(y-\bar{y}),\\
	\Sigma_{xx}' & = \Sigma_{xx} - \Sigma_{xy}\Sigma_{yy}^{-1}\Sigma_{yx},
  \end{split}
\end{equation}
respectively
\cite{anderson1979optimal}.

If both $\mathbf{x}_{n+1|n}$ and $\phi^n(\mathbf{x}_{n+1|n})$ are Gaussian,
one can apply
Eq.~(\ref{eq:conditionedmoments})
with
$X=\mathbf{x}_{n+1|n}$,
$Y = \phi^n(\mathbf{x}_{n+1|n})+\eta_{n+1}$
and $y = \mathbf{y}_{n+1}$
to obtain 
the first two moments of $\mathbf{x}_{n+1|n} \vert \mathbf{y}_{n+1} = \mathbf{x}_{n+1|n+1}$.
However
$\phi^n(\mathbf{x}_{n+1})$ is 
not a Gaussian,
unless $\phi^n$ is a linear function and
$\mathbf{x}_{n+1}$ is Gaussian.
As in the case of the time update,
we consider two Gaussian approximations:

\subsubsection{Linear Gaussian approximation}
The Taylor approximation of
\begin{equation*}
  \phi^n(\mathbf{x}_{n+1|n}) \simeq 
  \phi^n(\bar{\mathbf{x}}_{n+1|n})
  +\nabla \phi^n\vert_{ \bar{\mathbf{x}}_{n+1|n}}(\mathbf{x}_{n+1|n}-\bar{\mathbf{x}}_{n+1|n}),
\end{equation*}
which is Gaussian,
is used in place of
$\phi^n(\mathbf{x}_{n+1|n})$.
In this case, we use Eq.~(\ref{eq:conditionedmoments}) to obtain
\begin{equation}
  \begin{split}
	\label{eq:LGA}
\bar{\mathbf{x}}_{n+1|n+1}
& \simeq \bar{\mathbf{x}}_{n+1|n} + {G}_\mathbf{x}
 \left(\mathbf{y}_{n+1} - \phi^n(\bar{\mathbf{x}}_{n+1|n}) \right), \\
\mathbf{C}_{n+1|n+1}  
& \simeq \mathbf{C}_{n+1|n} -  
{G}_\mathbf{x} \nabla \phi^n\vert_{ \bar{\mathbf{x}}_{n+1|n}} \mathbf{C}_{n+1|n}
  \end{split}
\end{equation}
where
\begin{equation*}
  \begin{split}
&	{G}_\mathbf{x} \equiv
\mathbf{C}_{n+1|n}
(\nabla \phi^n\vert_{ \bar{\mathbf{x}}_{n+1|n}})^T  \\
& \qquad \left( 
\nabla \phi^n\vert_{ \bar{\mathbf{x}}_{n+1|n}}
\mathbf{C}_{n+1|n}
(\nabla \phi^n\vert_{ \bar{\mathbf{x}}_{n+1|n}})^T + R_{n+1} \right)^{-1}.
  \end{split}
\end{equation*}

\subsubsection{Point-based Gaussian approximation}
Let $\sum_j \lambda_j \delta_{ \mathbf{x}_{n+1|n}^j}$ be a discrete measure distributed according to
the probability distribution of $\mathbf{x}_{n+1|n}$
then $\sum_j \lambda_j \delta_{\phi^n(\mathbf{x}_{n+1|n}^j)}$
is distributed according to
the probability distribution of $\phi^n(\mathbf{x}_{n+1|n})$.
We approximate $\phi^n(\mathbf{x}_{n+1|n})$ by a Gaussian
whose mean and covariance are those obtained from
$\sum_j \lambda_j \delta_{ \phi^n(\mathbf{x}_{n+1|n}^j)}$.
In this case, we use Eq.~(\ref{eq:conditionedmoments}) to obtain
\begin{equation}
  \begin{split}
	\label{eq:PGA}
\bar{\mathbf{x}}_{n+1|n+1} & \simeq \bar{\mathbf{x}}_{n+1|n}  
+ {L}_\mathbf{x}\left(\mathbf{y}_{n+1} - \mathbf{z} \right), \\
\mathbf{C}_{n+1|n+1} & \simeq \mathbf{C}_{n+1|n} - {L}_\mathbf{x} P_{\mathbf{x}\mathbf{z}}^T
  \end{split}
\end{equation}
where
\begin{equation*}
  \begin{split}
	{L}_\mathbf{x} & \equiv P_{\mathbf{x}\mathbf{z}}\left(P_{\mathbf{z}\mathbf{z}}+R_{n+1}\right)^{-1}, \\
	\mathbf{z} & \equiv  \sum_j \lambda_j { \phi^n(\mathbf{x}_{n+1|n}^j)}, \\ 
	P_{\mathbf{x}\mathbf{z}} & \equiv \sum_j \lambda_j 
	\Big( {\mathbf{x}^j_{n+1|n}} -  \sum_j \lambda_j { \mathbf{x}_{n+1|n}^j}\Big)
	\left( \phi^n({\mathbf{x}^j_{n+1|n}}) - \mathbf{z} \right)^T, \\
	P_{\mathbf{z}\mathbf{z}} & \equiv \sum_j \lambda_j 
	\left( \phi^n({\mathbf{x}^j_{n+1|n}}) - \mathbf{z} \right)
	  \left( \phi^n({\mathbf{x}^j_{n+1|n}}) - \mathbf{z} \right)^T.
  \end{split}
\end{equation*}

In addition to the above two methods,
we consider an algorithm that does not use point-based approximation,
and does not require a Gaussian assumption regarding
$\phi^n(\mathbf{x}_{n+1})$.
It is motivated by 
the variational data assimilation
widely used in weather forecasting 
\cite{Fisher05}.

\subsubsection{Variational Gaussian approximation} 
Let the probability density function of a centered Gaussian with covariance $R_{n+1}$, be denoted by  $\Theta(\mathbf{y}_{n+1},R_{n+1})$.
Eq.~(\ref{eq:bayesrule}) then implies that
$\mathbb{P}(X|Y)  = {\mathbb{P}(X)\mathbb{P}(Y|X) }/{\mathbb{P}(Y) }$
and
\begin{equation*}
  \begin{split}
  \mathbb{P}(\mathbf{x}_{n+1|n+1}) & \propto \mathbb{P}(\mathbf{x}_{n+1|n})\, 
  \Theta\left(\mathbf{y}_{n+1}-\phi^n(\mathbf{x}_{n+1}), R_{n+1} \right) \\
& \propto \exp(-\mathbf{J}_{n+1|n+1}(\mathbf{x_{n+1}}))
  \end{split}
\end{equation*}
where the misfit function $ \mathbf{J}$ is given by
\begin{equation}
  \begin{split}
	\label{eq:misfit}
	& \mathbf{J}_{n+1|n+1}(\mathbf{x}_{n+1})\\ 
	& \!\!\! \equiv \frac{1}{2}\left \lbrace
	\parallel \mathbf{x}_{n+1}-\bar{\mathbf{x}}_{n+1|n}\parallel_{\mathbf{C}_{n+1|n}}^2
	+ \parallel \mathbf{y}_{n+1} - \phi^n (\mathbf{x}_{n+1})\parallel_{R_{n+1}}^2 \right\rbrace.
  \end{split}
\end{equation}
Here the notation $\parallel X \parallel_\Sigma^2 \equiv X^T \Sigma^{-1} X$ 
is used for a positive definite quadratic form with matrix $\Sigma$.
Therefore a Gaussian approximation of $\mathbf{x}_{n+1|n+1}$ is equivalent to 
making a quadratic approximation of
\begin{equation*}
	\mathbf{J}_{n+1|n+1}(\mathbf{x}_{n+1})
	 \simeq \frac{1}{2} \parallel \mathbf{x}_{n+1}-\bar{\mathbf{x}}_{n+1|n+1} 
	\parallel_{\mathbf{C}_{n+1|n+1}}^2
	+\text{const}.
\end{equation*}
The variational method approximates
$\bar{\mathbf{x}}_{n+1|n+1}$ 
by the minimizer of 
$\mathbf{J}_{n+1|n+1}$
and 
$\mathbf{C}_{n+1|n+1}$
by the inverse of 
the Hessian of 
the misfit function 
at $\bar{\mathbf{x}}_{n+1|n+1}$, i.e.,
\begin{equation}
  \begin{split}
	\label{eq:VDA}
\bar{\mathbf{x}}_{n+1|n+1}
& \simeq \text{minimizer of Eq.~(\ref{eq:misfit})}, \\
\mathbf{C}_{n+1|n+1}
& \simeq \left(\nabla \nabla \mathbf{J}_{n+1|n+1}\vert_{\bar{\mathbf{x}}_{n+1|n+1}} \right)^{-1}.
  \end{split}
\end{equation}

\subsection{Construction of Conventional Gaussian Filters}
\label{subsec:CGF}
We can choose one 
from
the two UQ methods 
(Eqs.~(\ref{eq:momentmapping}), (\ref{eq:momentmapping2}) with $n'=n$)
and independently one 
from
the three DA methods
(Eqs.~(\ref{eq:LGA}), (\ref{eq:PGA}), (\ref{eq:VDA}))
to construct a conventional filter.
In this paper we intend to make the UQ and DA methods consistent, if possible,
and not to simultaneously use
the non-point-based algorithm and point-based algorithm.
As a result,
we define
the linear Gaussian filter (\textbf{LGF})
as the combination of UQ with a linear Gaussian approximation 
and DA with a linear Gaussian approximation;
the variational Gaussian filter (\textbf{VGF})
as the combination of UQ with linear Gaussian approximation 
and DA with a variational Gaussian approximation;
the cubature Gaussian filter (\textbf{CGF})
and the particle Gaussian filter (\textbf{PGF})
as the combination of UQ with a point-based Gaussian approximation 
and DA with a point-based Gaussian approximation, 
for which cubature measure and empirical measure are employed respectively.

\section{Gaussian Smoothing Filters}
\label{sec:variationalsmoothing}
\noindent
By combining the 
DA methods presented in subsection~\ref{subsec:SMU}
and the 
UQ methods presented in subsection~\ref{subsec:STU},
we develop a number of Gaussian approximation smoothing filters
in subsection~\ref{subsec:GSF}.

\subsection{Measurement Update ($\mathcal{X}_{n|n}\Rightarrow \mathcal{X}_{n|n+1}$)}
\label{subsec:SMU}
\noindent
The methodology for the measurement update in a Gaussian smoothing filter is the same 
as in the case of conventional Gaussian filtering,
except for the use of  $\mathcal{X}_{n|n'}$ in place of $\mathbf{x}_{n+1|n'}$
(hence $\bar{\mathcal{X}}_{n|n'}$ and $\mathcal{C}_{n|n'}$ in place of 
$\bar{\mathbf{x}}_{n+1|n'}$ and $\mathbf{C}_{n+1|n'}$ respectively)
and $\Psi^n \equiv \phi^n \circ \Phi^n$ 
(due to $\mathbf{y}_{n+1} = \Psi^n(\mathcal{X}_n)+\eta_{n+1}$)
in place of $\phi^n$.
We mention that 
$\Psi^n$
in smoothing filter
might be a nonlinear function even when $\phi^n$ is linear.

\subsubsection{Linear Gaussian approximation}
As with Eq.~(\ref{eq:LGA}), we obtain
\begin{equation}
  \begin{split}
	\label{eq:SLGA}
\bar{\mathcal{X}}_{n|n+1}
& \simeq \bar{\mathcal{X}}_{n|n} + {G}_\mathcal{X}
 \left(\mathbf{y}_{n+1} - \Psi^n(\bar{\mathcal{X}}_{n|n}) \right), \\
\mathcal{C}_{n|n+1}  
& \simeq \mathcal{C}_{n|n} -  
{G}_\mathcal{X} \nabla \Psi^n\vert_{ \bar{\mathcal{X}}_{n|n}} \mathcal{C}_{n|n}
  \end{split}
\end{equation}
where
\begin{equation*}
  \begin{split}
&	{G}_\mathcal{X} \equiv
\mathcal{C}_{n|n}
(\nabla \Psi^n\vert_{ \bar{\mathcal{X}}_{n|n}})^T  \\
& \qquad \left( 
\nabla \Psi^n\vert_{ \bar{\mathcal{X}}_{n|n}}
\mathcal{C}_{n|n}
(\nabla \Psi^n\vert_{ \bar{\mathcal{X}}_{n|n}})^T + R_{n+1} \right)^{-1}.
  \end{split}
\end{equation*}

\subsubsection{Point-based Gaussian approximation}
Recall that
$\sum_{j} \Lambda_j \delta_{\mathcal{X}^j_{n|n'}}$
denotes a discrete measure distributed according to the probability distribution of 
$\mathcal{X}_{n|n'}$.
As with Eq.~(\ref{eq:PGA}), we obtain
\begin{equation}
  \begin{split}
	\label{eq:SPGA}
\bar{\mathcal{X}}_{n|n+1} & \simeq \bar{\mathcal{X}}_{n|n}  
+ {L}_\mathcal{X}\left(\mathbf{y}_{n+1} - \mathcal{Z} \right), \\
\mathcal{C}_{n|n+1} & \simeq \mathcal{C}_{n|n} - {L}_\mathcal{X} P_{\mathcal{X}\mathcal{Z}}^T
  \end{split}
\end{equation}
where
\begin{equation*}
  \begin{split}
	{L}_\mathcal{X} & \equiv P_{\mathcal{X}\mathcal{Z}}
	\left(P_{\mathcal{Z}\mathcal{Z}}+R_{n+1}\right)^{-1}, \\
	\mathcal{Z} & \equiv  \sum_j \Lambda_j { \Psi^n(\mathcal{X}_{n|n}^j)}, \\ 
	P_{\mathcal{X}\mathcal{Z}} & \equiv \sum_j \Lambda_j 
	\Big( {\mathcal{X}^j_{n|n}} -  \sum_j \Lambda_j { \mathcal{X}_{n|n}^j}\Big)
	\left( \Psi^n({\mathcal{X}^j_{n|n}}) - \mathcal{Z} \right)^T, \\
	P_{\mathcal{Z}\mathcal{Z}} & \equiv \sum_j \Lambda_j 
	\left( \Psi^n({\mathcal{X}^j_{n|n}}) - \mathcal{Z} \right)
	  \left( \Psi^n({\mathcal{X}^j_{n|n}}) - \mathcal{Z} \right)^T.
  \end{split}
\end{equation*}

\subsubsection{Variational Gaussian approximation}
Let 
\begin{equation}
  \begin{split}
	\label{eq:misfit2}
	& \mathcal{J}_{n|n+1}(\mathcal{X}_n) \\
	& = \frac{1}{2}\left\lbrace
	\parallel \mathcal{X}_n-\bar{\mathcal{X}}_{n|n}\parallel_{\mathcal{C}_{n|n}}^2
	+ \parallel \mathbf{y}_{n+1} 
	- \Psi^n(\mathcal{X}_n)\parallel_{R_{n+1}}^2 \right\rbrace
  \end{split}
\end{equation}
be the misfit function.
As with Eq.~(\ref{eq:VDA}), we obtain
\begin{equation}
  \begin{split}
\label{eq:SVDA}
\bar{\mathcal{X}}_{n|n+1}
& \simeq \text{minimizer of Eq.~(\ref{eq:misfit2})}, \\
\mathcal{C}_{n|n+1}
& \simeq \left(\nabla \nabla \mathcal{J}_{n|n+1}\vert_{\bar{\mathcal{X}}_{n|n+1}} \right)^{-1}.
  \end{split}
\end{equation}

\subsection{Time Update ($\mathcal{X}_{n|n+1} \to \mathbf{x}_{n+1|n+1}$)}
\label{subsec:STU}
We apply Eq.~(\ref{eq:momentmapping}) or 
Eq.~(\ref{eq:momentmapping2})
with $n'=n+1$
to perform the UQ of the smoothing filter.
Unlike the case of $n'=n$ in the conventional filtering,
$\mathbb{E}({\xi}_{n|n+1})\neq \mathbf{0}$ 
and $\Sigma_{n|n+1}$ is not block diagonal.

\subsection{Construction of Gaussian Smoothing Filters}
\label{subsec:GSF}
We can choose 
one 
from
the three DA methods
(Eqs.~(\ref{eq:SLGA}), (\ref{eq:SPGA}), (\ref{eq:SVDA}))
and independently
one 
from
the two UQ methods 
(Eqs.~(\ref{eq:momentmapping}), (\ref{eq:momentmapping2}) with $n'=n+1$)
to construct a conventional filter.
In this paper we intend to develop a smoothing filter 
modelled upon a given conventional filter
or to make a one-to-one correspondence between conventional filters and smoothing filters.
As a result,
we define
the linear Gaussian smoothing filter (\textbf{LGSF})
as the combination of 
DA with a linear Gaussian approximation and UQ with a linear Gaussian approximation;
the variational Gaussian smoothing filter (\textbf{VGSF})
as the combination of 
DA with a variational Gaussian approximation and
UQ with a linear Gaussian approximation;
the cubature Gaussian smoothing filter (\textbf{CGSF})
and the particle Gaussian smoothing filter (\textbf{PGSF})
as the combination of DA with a point-based Gaussian approximation 
and UQ with a point-based Gaussian approximation, 
for which cubature measure and empirical measure are employed respectively.
We mention that LGSF, VGSF, CGSF and PGSF
correspond to LGF, VGF, CGF and PGF, respectively, and vice versa.
We also mention that the computational difference between two corresponding filters lies
at the DA step.

\section{Practical Implementation}
\label{sec:comparison}
In this section some practical issues,  encountered in implementing smoothing filters,
are resolved.

\begin{enumerate}
  \item
In the case that the forward model derives from an approximation to the stochastic differential equation,
\begin{equation}
  \label{eq:sde}
  d\mathbf{x}(t) = b(t,\mathbf{x}(t))dt+s(t,\mathbf{x}(t))\,dB(t),
\end{equation}
where
$b \in \mathbb{R}^{d}$ is the drift,
$s \in \mathbb{R}^{d \times N}$ is the volatility and 
$B=(B_1,\cdots,B_{N})$ is the set of independent Brownian motions,
describing the evolution of the 
underlying system to be estimated,
one constructs
the forward model 
of 
Eq.~(\ref{eq:fmodel}) 
as follows.
Let $\delta t > 0$ be the numerical simulation time step
and let 
the observations arrive at the times
$\Delta t = M\times \delta t$.
A finite difference approximation
of Eq.~(\ref{eq:sde})
using the Euler-Maruyama or Milstein method 
\cite{kloeden2011numerical}
yields
\begin{equation}
  \label{eq:discretesde}
  \mathbf{x}_{n,m+1}  = \mathcal{F}^{n,m}(\mathbf{x}_{n,m}, w_{n,m}), 
  \quad w_{n,m} \sim \mathcal{N}( \mathbf{0},Q_{n,m} )
\end{equation}
where
$\mathbf{x}_{n,m}$ denotes an approximation of $\mathbf{x}(n\Delta t + m\delta t)$.
The repeated application of Eq.~(\ref{eq:discretesde})
from $m=0$ to $m=M-1$
defines 
$\Phi^n(\cdot)$ of the forward model,
the mapping from
$\mathbf{x}_n= \mathbf{x}_{n,0}$
to
$\mathbf{x}_{n+1}= \mathbf{x}_{n,M}$,
along with 
the augmented vector
$ \xi_n = 
\left[
\begin{array}{c}
  w_{n,0}  \\ 
 \cdots \\
  w_{n,M-1}
\end{array}\right] $
and the block diagonal matrix
$ \Gamma_n 
=\left[
\begin{array}{ccc}
  Q_{n,0} &  & \mathbf{0} \\ 
   & \cdots &  \\ 
  \mathbf{0} &  & Q_{n,M-1}
\end{array}\right]$.

\item
Let $\sum_j \omega_j \delta_{x^j}$ be 
a cubature measure
or
an empirical measure
approximating the standard Gaussian.
Then $\sum_j \omega_j \delta_{ m + S x^j}$, for which $S$ satisfies $C = SS^T$,
becomes an approximation for $\mathcal{N}(m,C)$.

\item
Some cubature formulae with respect to the standard Gaussian 
can be found in
\cite{cools1993monomial,
victoir2004asymmetric,
xiu2008numerical}.
In one dimension, 
a cubature measure is more commonly referred to as a quadrature measure.
A general multi-dimensional Gaussian cubature can be constructed
via the tensor product of a quadrature formula
\cite{arasaratnam2009cubature}.
Using Gauss-Hermite quadrature with support size $s = (r+1)/2$ for degree $r$,
one can develop
a cubature formula of degree $r$ with respect to a $k$-dimensional standard Gaussian
whose support size is $s^k$.
Because the computational cost increases 
as 
the support size of the discrete measure increases,
it is important to use cubature measure supported on a smaller set. 
For the numerical simulations performed in the next section
we use 
the standard Gaussian cubature formula of degree $3$ and $5$
introduced in \cite{jia2012high},
whose support size is $2k$ and $2k^2+1$, respectively.

\item
The Broyden-Fletcher-Goldfarb-Shanno (BFGS) iterative method 
\cite{mordecai2003nonlinear},  \cite{matlab}, is used
to solve the nonlinear optimisation problem for the variational Gaussian approximation. Numerical derivatives are employed as the dimension of our test cases is low. More generally an adjoint derivative method could be used for greater efficiency.

\end{enumerate}

\section{Numerical Simulations}
\label{eq:numericalsimulations}
In this section the feasibility of our proposed smoothing filters  is investigated.
The performance of 
LGSF, VGSF, CGSF and PGSF is compared with  that of LGF, VGF, CGF and PGF.
The metric used to compare the performance of various
filters is the root mean square error
(RMSE).
The 
RMSE
between $A=\{ A_i\}_{i=1}^N$ and $B=\{ B_i\}_{i=1}^N$
is defined by
\begin{equation*}
  \textstyle
  \text{RMSE}(A,B) = \sqrt{ \frac{1}{N} \sum_{i=1}^N \vert A_i - B_i \vert^2 }
\end{equation*}
where $A_i$ and $B_i$ are vectors.
The index $i$ will specify either simulation number or time step.

The following comparison study uses various examples 
with different starting states and problem data.
In the examples we studied, 
we see the smoothing filters generally yield more accurate estimations 
than the corresponding conventional filters.
The test examples consist of 
(i) a bistable system 
(subsection~\ref{subsec:bistable}),
(ii) a prototypical chaotic dynamical system
(subsection~\ref{subsec:L63})
and (iii) a target tracking problem 
(subsection~\ref{subsec:target}).

\subsection{Bistable System}
\label{subsec:bistable}
We consider the one-dimensional differential equation 
\begin{equation}
  \begin{split}
	\label{eq:1dsde}
	d\mathbf{x}(t) = \beta \mathbf{x}(1-\mathbf{x}^2) dt +\sigma dB(t), \quad \beta > 0.
  \end{split}
\end{equation}
The deterministic equation 
with $\sigma=0$
has two stable equilibria, $-1$ and $+1$, and one unstable
equilibrium, $0$.
In the deterministic case the process $\mathbf{x}(t)$ is distributed around one of the stable equilibria.
The stochastic system ($\sigma \neq 0$) 
however  
shows sudden transitions between the two stable equilibria
due to the presence of
random perturbation
\cite{freidlin2012random}.
The Euler approximation of Eq.~(\ref{eq:1dsde})
is used to produce
\begin{equation}
  \label{eq:dissde}
  \mathbf{x}_{n,m+1} = 
  \mathbf{x}_{n,m}+\delta t \times \beta \mathbf{x}_{n,m}(1-\mathbf{x}_{n,m}^2) 
  +\mathcal{N}(0,\sigma^2 \delta t)
\end{equation}
which corresponds to Eq.~(\ref{eq:discretesde}).

\begin{figure}
  \subfigure[One instance of state estimation of $\mathbf{x}_n$]{
\includegraphics[width=0.47\textwidth]{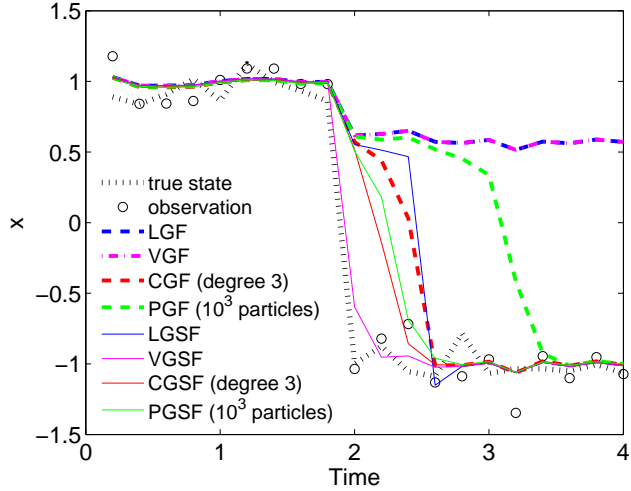}
\label{fig:dwp2a} 
}
\subfigure[The average of $50$ independent state estimations]{
\includegraphics[width=0.47\textwidth]{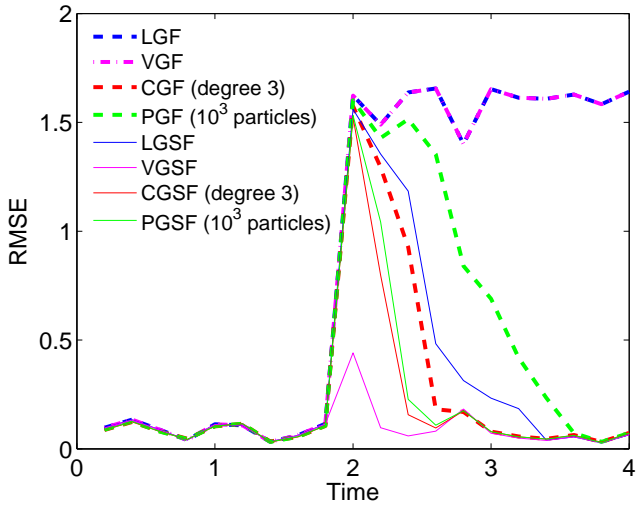}
\label{fig:dwp2b} 
}
\caption{
  The performance of various filters applied to the bistable system 
  with the identity measurement function.
}
\label{fig:dwp2} 
\end{figure}

We first study the case 
for which the measurement function 
is the identity function, i.e.,
\begin{equation}
  \begin{split}
	\label{eq:linobs}
	\phi^n(\mathbf{x}) & = \mathbf{x}.
  \end{split}
\end{equation}
For the system parameters
$\beta = 10$ and $\sigma=0.5$,
a single realisation 
of Eq.~(\ref{eq:1dsde})
starting at
$\mathbf{x}(0) = 0.8$
is simulated 
using Eq.~(\ref{eq:dissde})
with $\delta t = 0.01$
and regarded as the true state.
We see
this trajectory has
a jump from $1$ to $-1$ at $t=2.0$
(see Fig.~\ref{fig:dwp2a}).
We perform $50$ independent numerical approximations to estimate the true state.
In each case
the observational data is generated with 
$R_n =0.03$
at every
$\Delta t = 20 \times \delta t = 0.2$.
We then apply the conventional filters
as well as the smoothing filters
with the initial probability distribution
$\mathbf{x}_0 \sim \mathcal{N}(0.8, 0.02)$.
Fig.~\ref{fig:dwp2a} shows a representative case among our numerical simulations.
There, we depict the conditioned mean of the filtering solutions
together with the true state.
We are particularly interested in state estimates obtained from the filters 
since the transition takes place,
i.e., for $t \geq 2.0$.
In this case
the non-point-based conventional filters (LGF and VGF)
completely lose 
the true state.
The point-based conventional filters (CFG and PGF) 
eventually catch the trajectory 
but 
after a number of assimilation time steps.
Finally we see that
the smoothing filters build an accurate reconstruction of the state evolution despite the jump.
This is 
clearly due to 
the automatic importance sampling characteristic
described in subsection~(\ref{subsec:Importance}).
We here notice that the one from VGSF 
very quickly follows
the true state after the sudden transition happens.
Fig.~\ref{fig:dwp2b} shows the average RMSEs of $50$ state estimates
for each time step.

\begin{figure}[t!]
  \subfigure[Frequent observation]{
\includegraphics[width=0.43\textwidth]{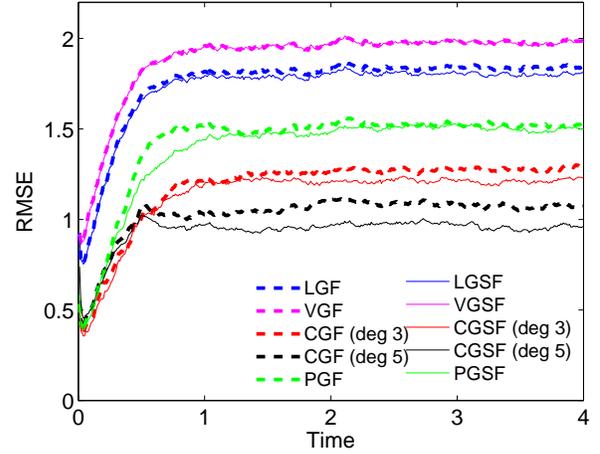}
\label{fig:dwp1a} 
}
  \subfigure[Sparse observation]{
\includegraphics[width=0.43\textwidth]{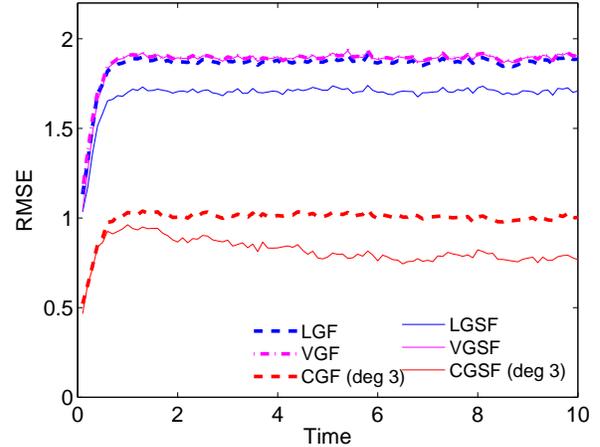}
\label{fig:dwp1b} 
}
\caption{
  The performance of various filters applied to the bistable system 
  with the shifted quadratic measurement function.
  For PGF and PGSF, empirical measures consisting of $10^3$ random samples are used.
}
\label{fig:dwp1} 
\end{figure}

\begin{figure}[t!]
  \subfigure[$x_1$]{
\includegraphics[width=0.43\textwidth]{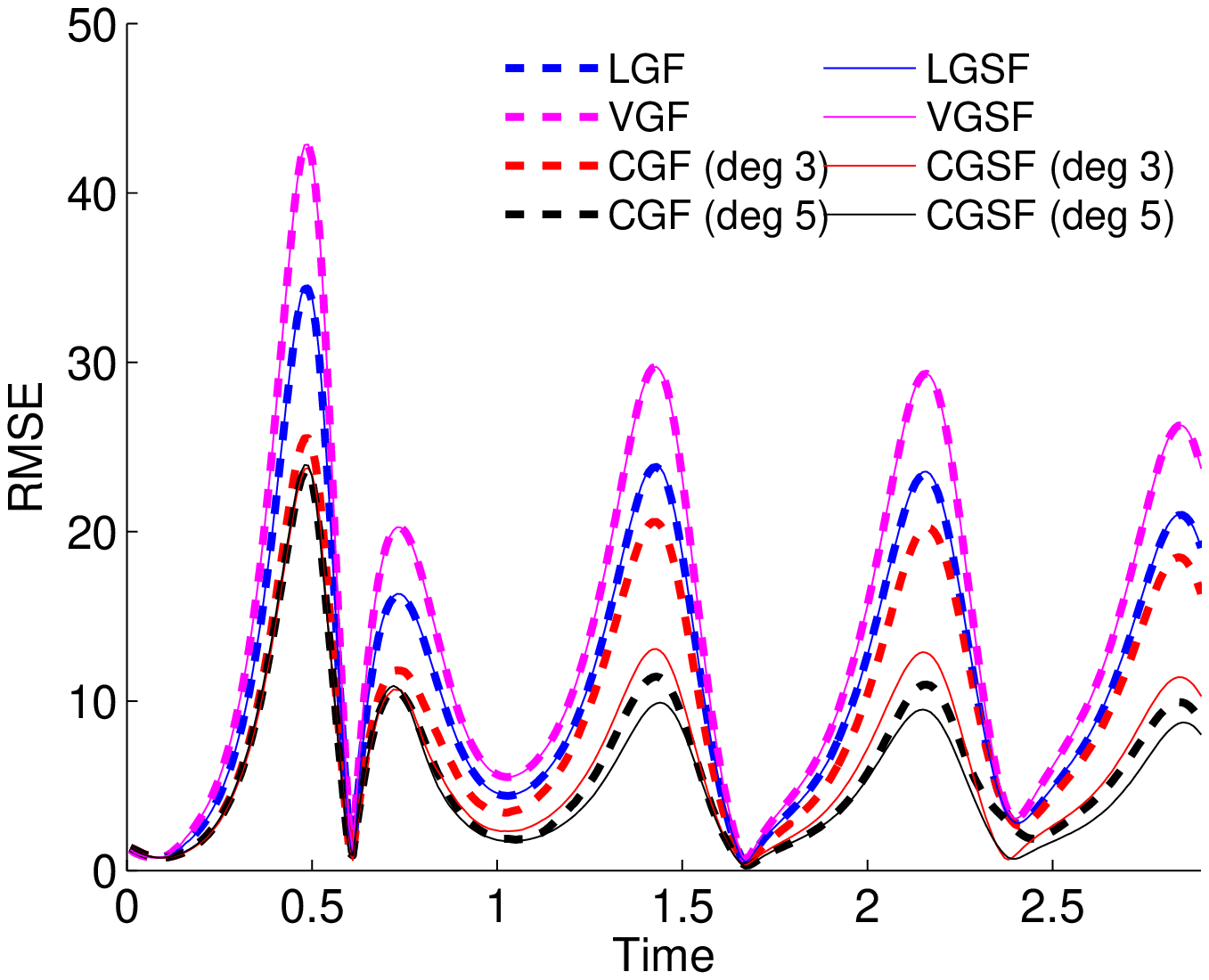} 
}
  \subfigure[$x_3$]{
\includegraphics[width=0.43\textwidth]{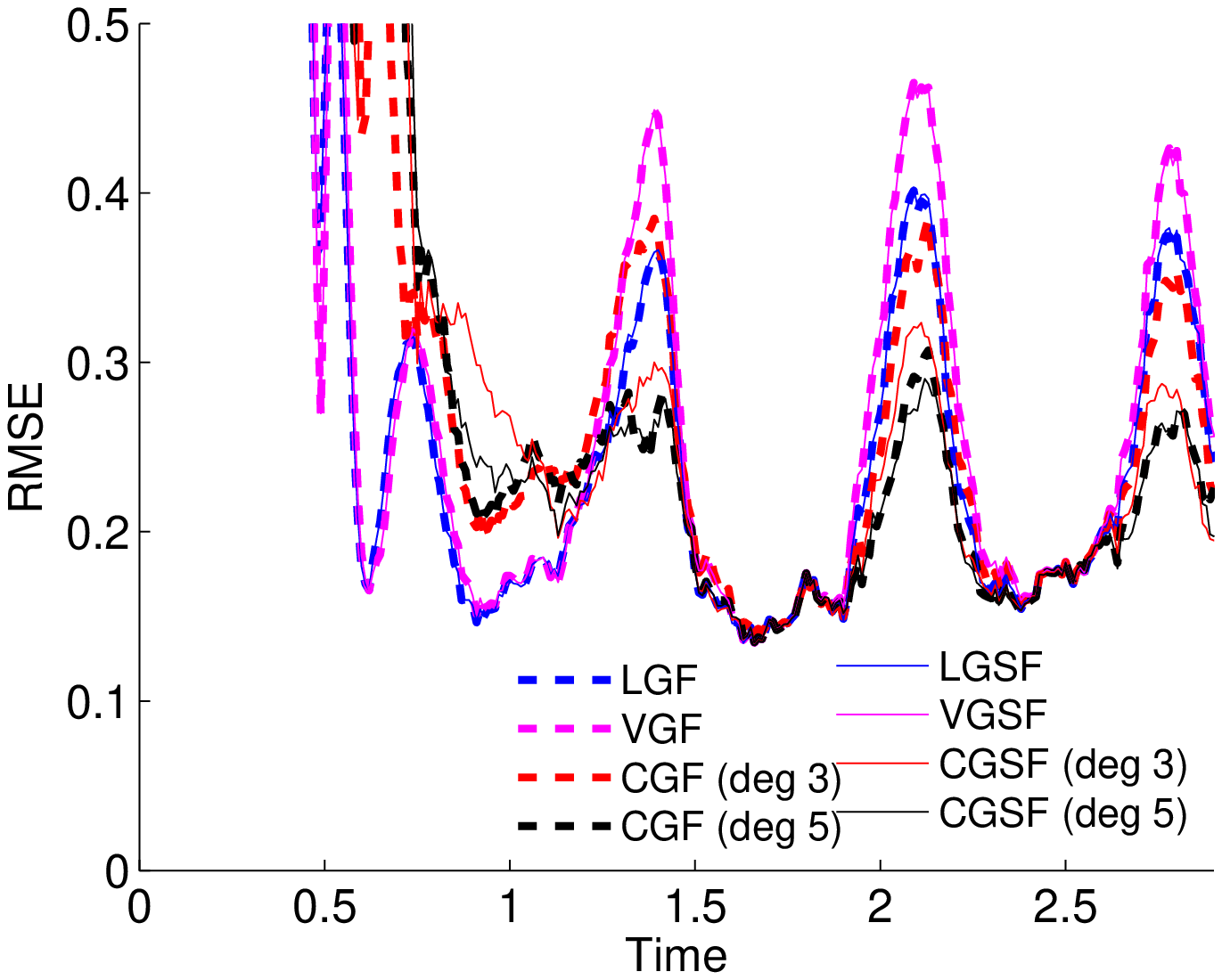} 
}
\caption{
  The performance of various filters applied to the Lorenz-$63$ model
  with the shifted quadratic measurement function.
  The plot for $x_2$ is very similar to that of $x_1$ and is not shown.
}
\label{fig:lorenz1} 
\end{figure}

We next study the case 
for which the measurement function 
is the square of the shifted distance from the origin
\begin{equation*}
  \phi^n(\mathbf{x})  = (\mathbf{x}-0.05)^2.
\end{equation*}
The observation
distinguishes the two stable equilibria marginally. 
Fig.~\ref{fig:dwp1} uses the system parameters $\beta = 5$,
$\sigma=0.5$,
the initial state $\mathbf{x}(0) = -0.2$,
the numerical simulation time step
$\delta t = 0.01$,
and
the observation noise covariance $R_n =1.0$.
Along with the initial condition $\mathbf{x}_0 \sim \mathcal{N}(0.8, 2.0)$,
the filters are applied at various inter-observation times $\Delta t = M \times \delta t$.
The average RMSEs
committed by each filter across $100$ independent simulations
are depicted
when $M=1$
(Fig.~\ref{fig:dwp1a})
and 
when $M=10$
(Fig.~\ref{fig:dwp1b}).
In this example one can see that 
the point-based conventional filters outperform 
the non-point-based conventional filters 
and that
the accuracy of smoothing filters are improved compared with
corresponding conventional filters.
Furthermore, as the time between two successive measurements increases, 
the smoothing filters become more accurate
compared with corresponding conventional filters.
This improvement of smoothing filters for temporally sparse observations
can also be understood from the importance sampling characteristic.

\subsection{Lorenz-$63$ System}
\label{subsec:L63}
Let $\mathbf{x}(t) = [x(t), y(t), z(t)]^T$
be the state vector.
We use the Euler approximation of the chaotic dynamical system
\begin{equation*}
  \begin{split}
  dx&=\sigma(y-x)dt+g_1 dB_1,\\
  dy&=(\rho x-y-xz)dt+g_2 dB_2,\\
  dz&=(xy-\beta z)dt+g_3 dB_3,
  \end{split}
\end{equation*}
with $\delta t = \Delta t = 0.01$
as the forward model 
\cite{lorenz1963deterministic,miller1999data}.
We choose the system parameters $\sigma=10$, $\rho=28$, $\beta = 8/3$ and $g_1=g_2=0$, $g_3=0.5$.
The starting state is
$\mathbf{x}(0) = [-0.2, \; -0.3, \; -0.5]^T$
and the initial condition is
$\mathbf{x}_0 \sim \mathcal{N}([1.35, \; -3, \; 6]^T, 0.35 I_3)$
where $I_3$ denotes the $3\times 3$ identity matrix.
The observation process is 
determined by the measurement function 
\begin{equation*}
  \begin{split}
	\phi^n(\mathbf{x}) & = \sqrt{ (x-0.5)^2+y^2+z^2 }
  \end{split}
\end{equation*}
and the noise covariance $R_n = 0.5$.
Fig.~\ref{fig:lorenz1} 
depicts the average RMSEs 
from $120$ simulations
for each component of system variables,
obtained from the conventional filters and smoothing filters.
The ordering of filtering accuracy 
among the different methods is very similar as the bistable system with frequent squared observation.

\subsection{Target Tracking}
\label{subsec:target}
Here we consider a model air-traffic monitoring scenario,
where an aircraft executes a maneuvering turn in a horizontal plane at an unknown turn rate 
$\Omega_n$ at time $n$.
The dynamical system is governed by the equation
\begin{equation*}
  \begin{split}
	& \mathbf{x}_{n+1} \\
& =
\left[
\begin{array}{ccccc}
  1 & \frac{\sin(\Omega_n\Delta t)}{\Omega_n} & 0 & \frac{\cos(\Omega_n\Delta t)-1}{\Omega_n} & 0 \\ 
  0 & \cos(\Omega_n\Delta t) & 0 & -\sin(\Omega_n\Delta t) & 0 \\ 
  0 & \frac{1-\cos(\Omega_n\Delta t)}{\Omega_n} & 1 & \frac{\sin(\Omega_n\Delta t)}{\Omega_n} & 0 \\ 
  0 & \sin(\Omega_n\Delta t) & 0 & \cos(\Omega_n\Delta t) & 0 \\ 
  0 & 0 & 0 & 0 & 1 
\end{array}\right]\mathbf{x}_n + \xi_n
  \end{split}
\end{equation*}
where $\mathbf{x}_n = [x_n, \; \dot{x}_n,\;  y_n,\;  \dot{y}_n,\;  \Omega_n]^T$;
$[x_n,\; y_n]$ and $[\dot{x}_n,\; \dot{y}_n]$ are the position and velocity 
of the aircraft
at time $n$;
$\Delta t$ is the time interval between two consecutive measurements;
the driving noise
$\xi_n \in \mathbb{R}^5$ is the zero mean Gaussian with covariance matrix
\begin{equation*}
 \Gamma_n = 
\left[\begin{array}{ccccc}
	\frac{\Delta t^3}{3} & \frac{\Delta t^2}{2} & 0 & 0 & 0 \\ 
	\frac{\Delta t^2}{2} & \Delta t & 0 & 0 & 0 \\ 
  0 & 0 & \frac{\Delta t^3}{3} & \frac{\Delta t^2}{2} & 0 \\ 
  0 & 0 & \frac{\Delta t^2}{2} & \Delta t & 0 \\ 
  0 & 0 & 0 & 0 & q \Delta t
\end{array}\right].
\end{equation*}
  Here the scalar parameter $q$ controls the random walk behavior of the turn rate
  from
$\Omega_{n+1} = \Omega_n + \mathcal{N}(0,q\Delta t)$.

We assume a radar is fixed at the origin of the plane and equipped to measure the range, $\rho_n$,
and bearing, $\theta_n$, at time $n$. Hence the observation process is
\begin{equation*}
\mathbf{y}_n =
\left[ \begin{array}{c}
  \rho_n \\
  \theta_n 
\end{array}\right]
+\eta_n
=
\left[
\begin{array}{c}
  \sqrt{x_n^2+y_n^2} \\
  \tan^{-1}\left(\frac{y_n}{x_n} \right)
\end{array}\right]
+\eta_n
\end{equation*}
where the measurement noise is $\eta_n \sim \mathcal{N}(\mathbf{0}, R_n)$
with 
\begin{equation*}
R_n =
\left[ \begin{array}{cc}
	\sigma_{\rho_n}^2 & 0 \\
	0 & \sigma_{\theta_n}^2
\end{array}\right].
\end{equation*}
Due to the inherent nonlinearity of the observation function,
target tracking is another problem 
suitable for testing
the performance of smoothing filters.

With the parameters
$\Delta t = 1$,
$q = 1.75\times 10^{-3}$,
$\sigma_{\rho_n}^2 = 10^2$,
$\sigma_{\theta_n}^2 = 10^{-5}$
and
\begin{equation*}
\mathbf{x}_0 \sim 
\mathcal{N}\left(
\left[ \begin{array}{c}
 10^3  \\
 3\times 10^2 \\
 10^3 \\
 0 \\
 -\frac{3\pi}{180}
\end{array}\right],
\left[ \begin{array}{ccccc}
 10^2 & 0 & 0 & 0 & 0  \\
 0 & 10 & 0 & 0 & 0  \\
 0 & 0 & 10^2 & 0 & 0  \\
 0 & 0 & 0 & 10 & 0  \\
 0 & 0 & 0 & 0 & 10^{-4}
\end{array}\right]
\right),
\end{equation*}
we perform 
$200$ independent simulations.
In each case 
the target trajectory,
whose initial state 
is 
an independent draw
from $\mathbf{x}_0$,
and the associated observations over 
$1 \leq n \leq 200$ time steps 
are randomly generated.
We then apply the filters to reconstruct the evolution of the dynamical variables.

\begin{figure}[t!]
\includegraphics[width=0.43\textwidth]{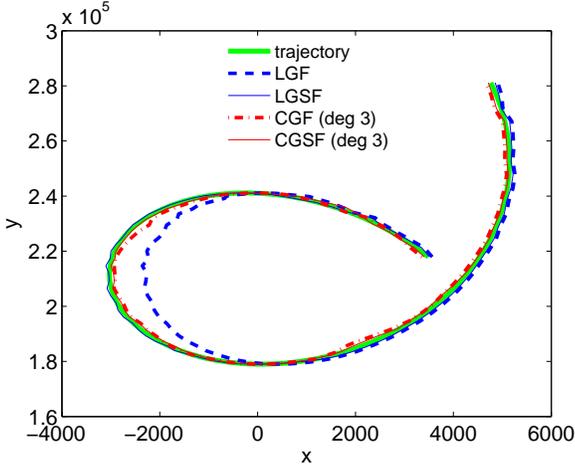}
\caption{
  A  trajectory of the aircraft together with 
  filtering estimates
  for $120 \leq n \leq 200$.
  The forward time  direction is counterclockwise.
}
\label{fig:bea1} 
\end{figure}

\begin{figure}[h!]
\subfigure[position]{
\includegraphics[width=0.43\textwidth]{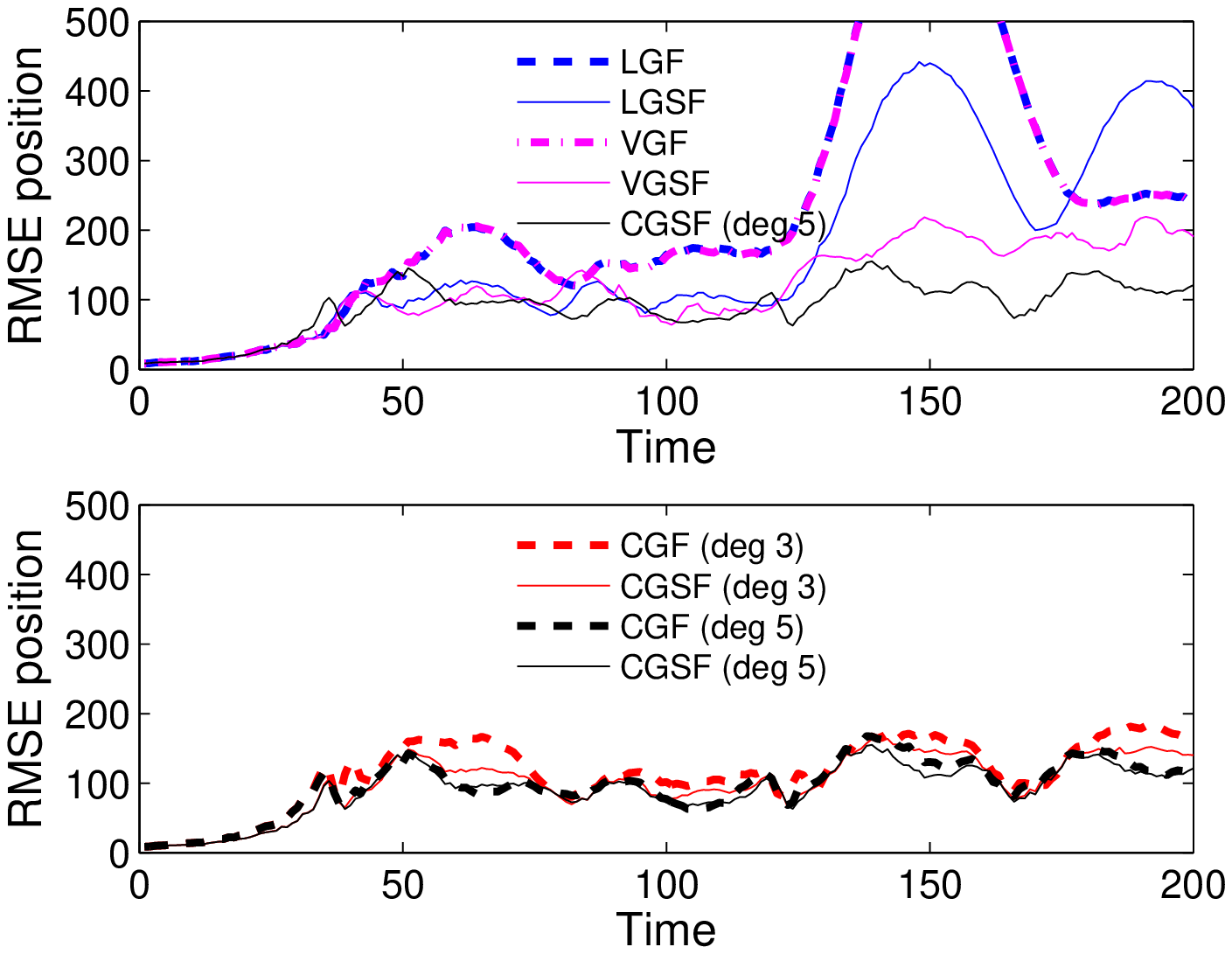}
\label{fig:bea2a} 
}
\subfigure[velocity]{
\includegraphics[width=0.43\textwidth]{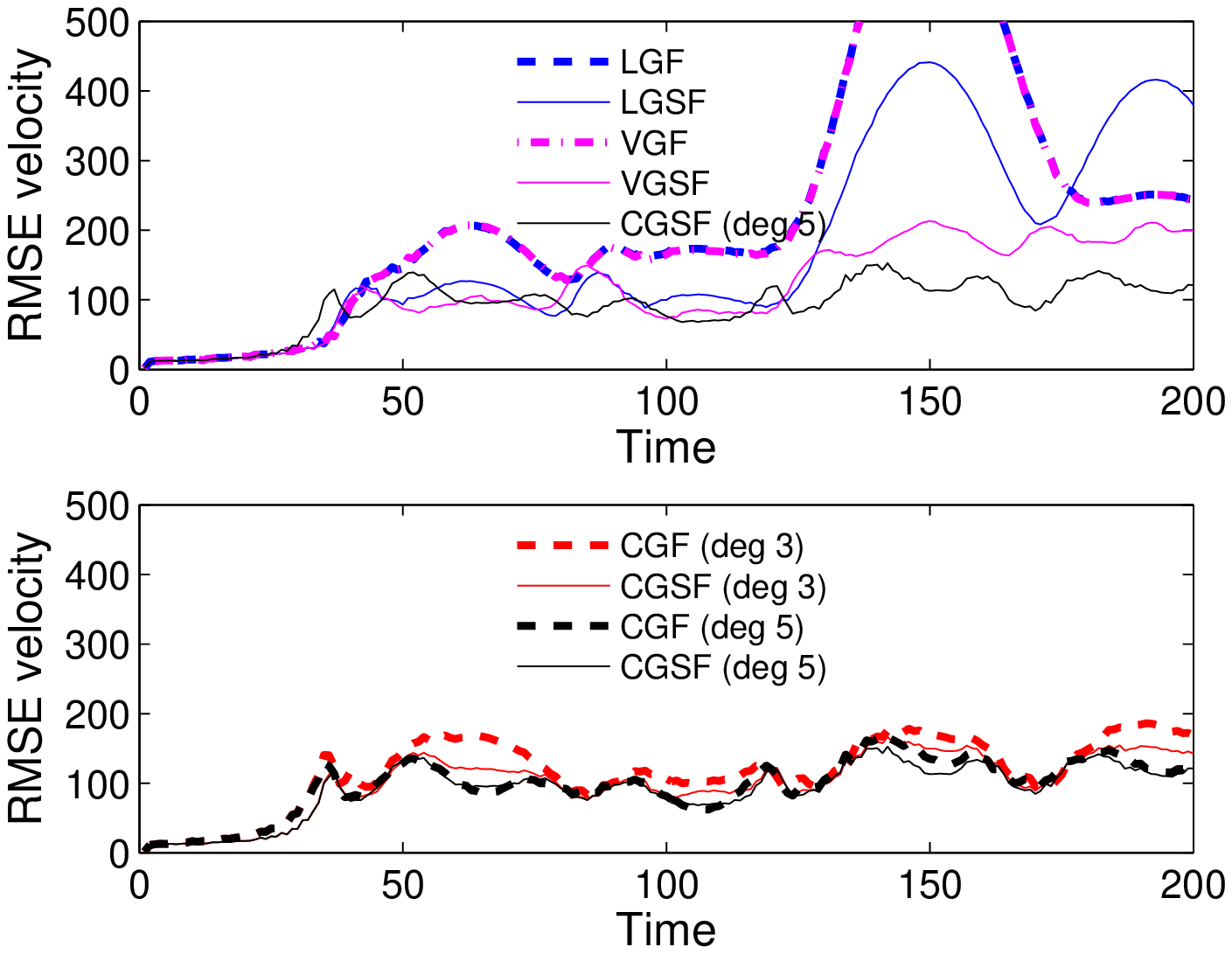} 
\label{fig:bea2b}
}
\subfigure[turn rate]{
\includegraphics[width=0.43\textwidth]{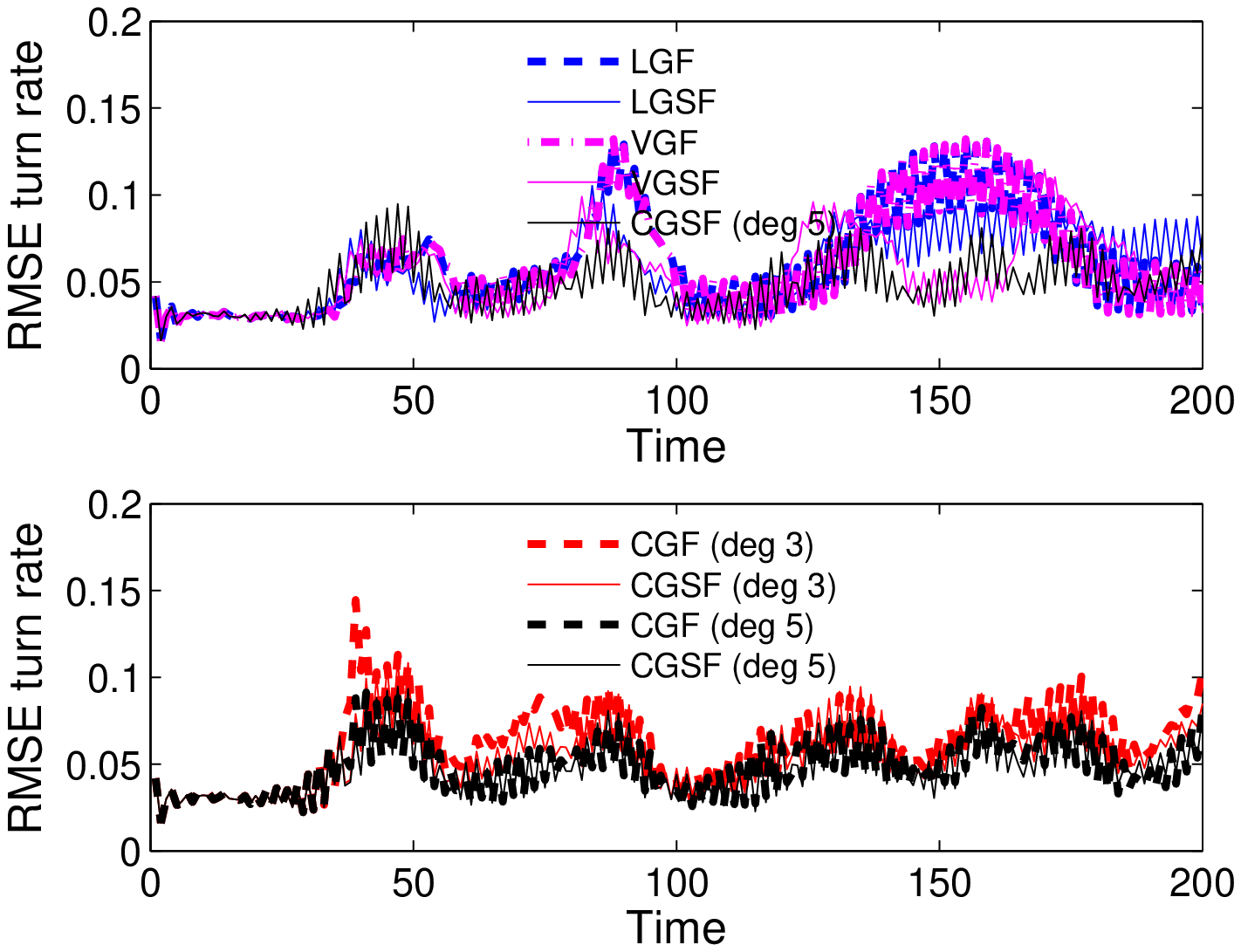} 
\label{fig:bea2c} 
}
\caption{
  The RMSEs between target and filtering estimates,
  obtained from averaging over $200$ independent simulations.
}
\label{fig:bea2} 
\end{figure}

\begin{figure*}[t!]
\begin{center}
\subfigure[position]{
\includegraphics[width=0.7\textwidth]{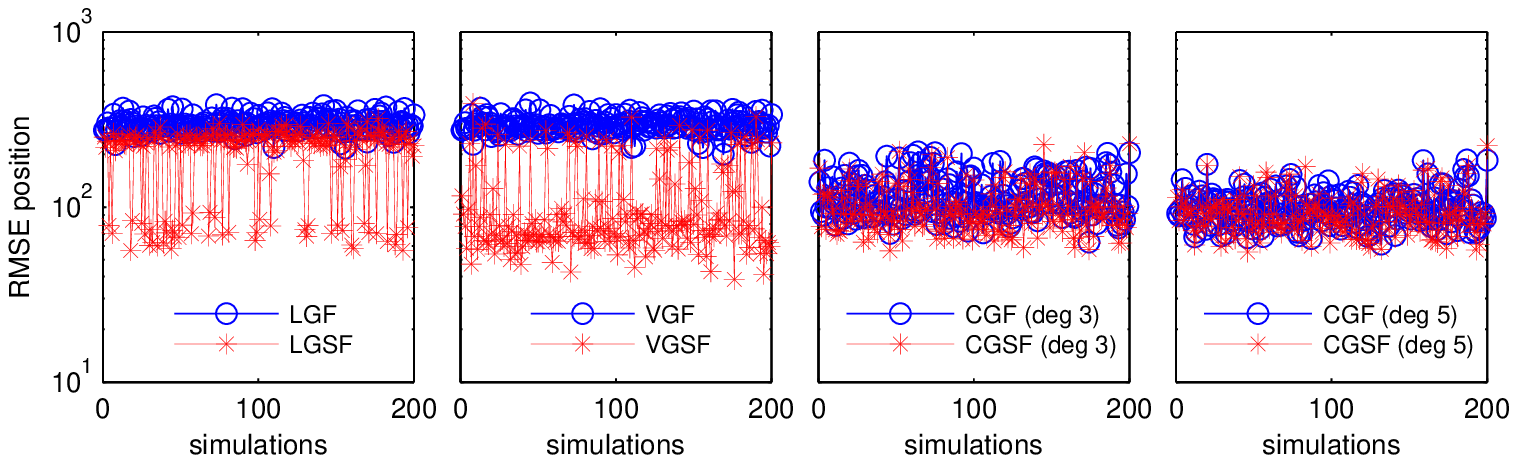}  
}
\subfigure[velocity]{
\includegraphics[width=0.7\textwidth]{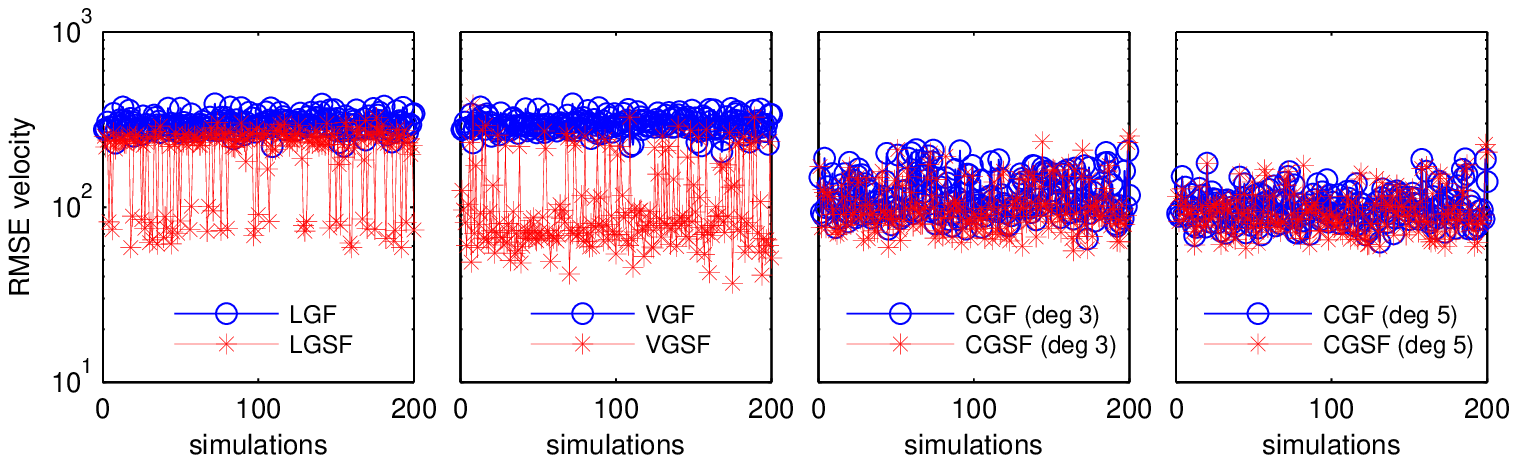}  
}
\subfigure[turn rate]{
\includegraphics[width=0.7\textwidth]{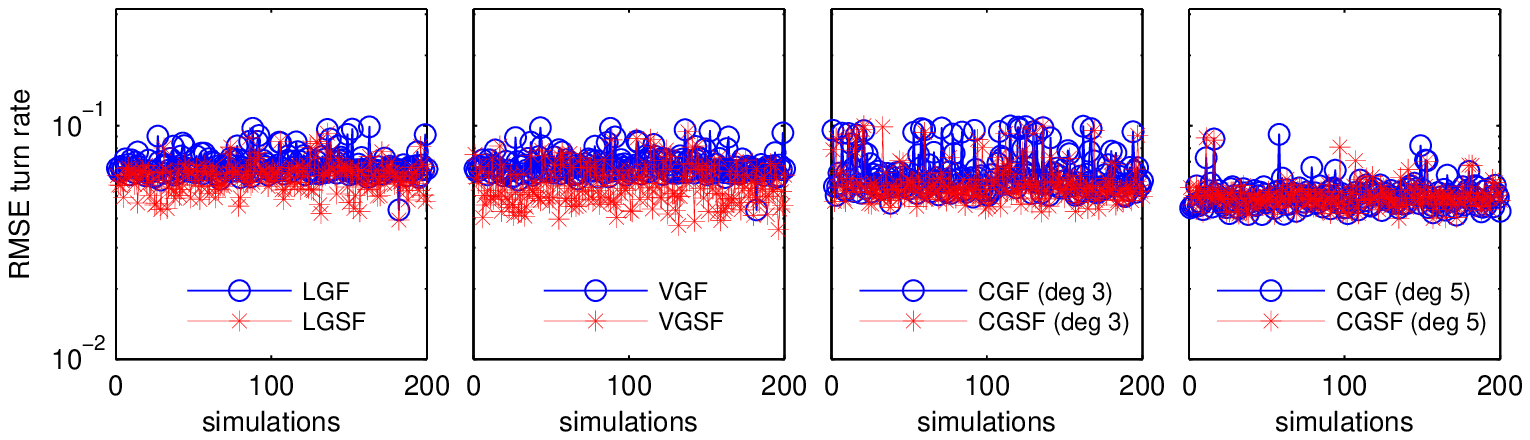} 
}
\end{center}
\caption{
  The RMSEs between target and filtering estimations,
  obtained from averaging over the time period $50 \leq n \leq 200$
  for each $200$ simulations.
}
\label{fig:bea3} 
\end{figure*}

\begin{figure}[t!]
\subfigure[position]
{\includegraphics[width=0.43\textwidth]{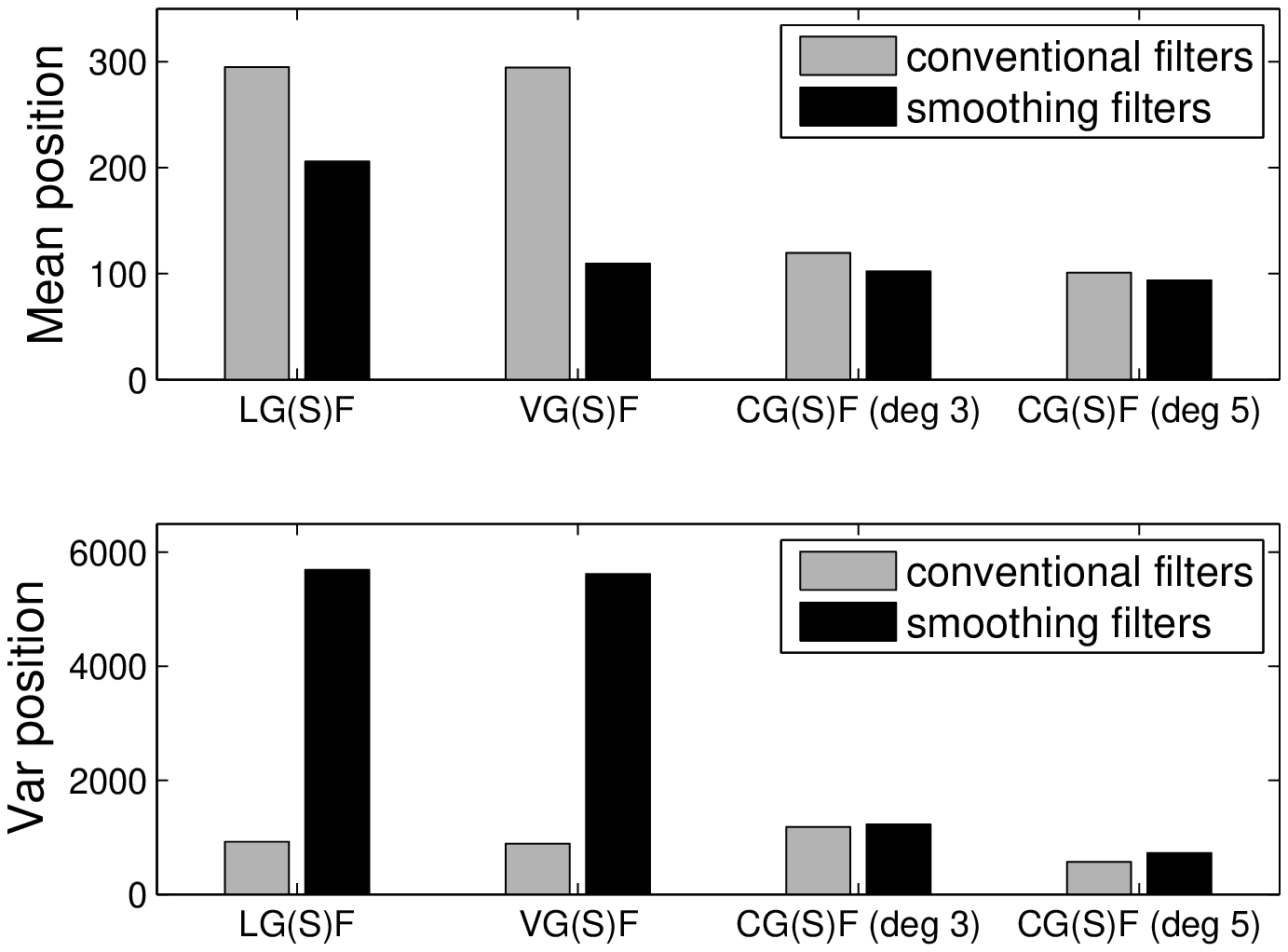} }
\subfigure[velocity]
{\includegraphics[width=0.43\textwidth]{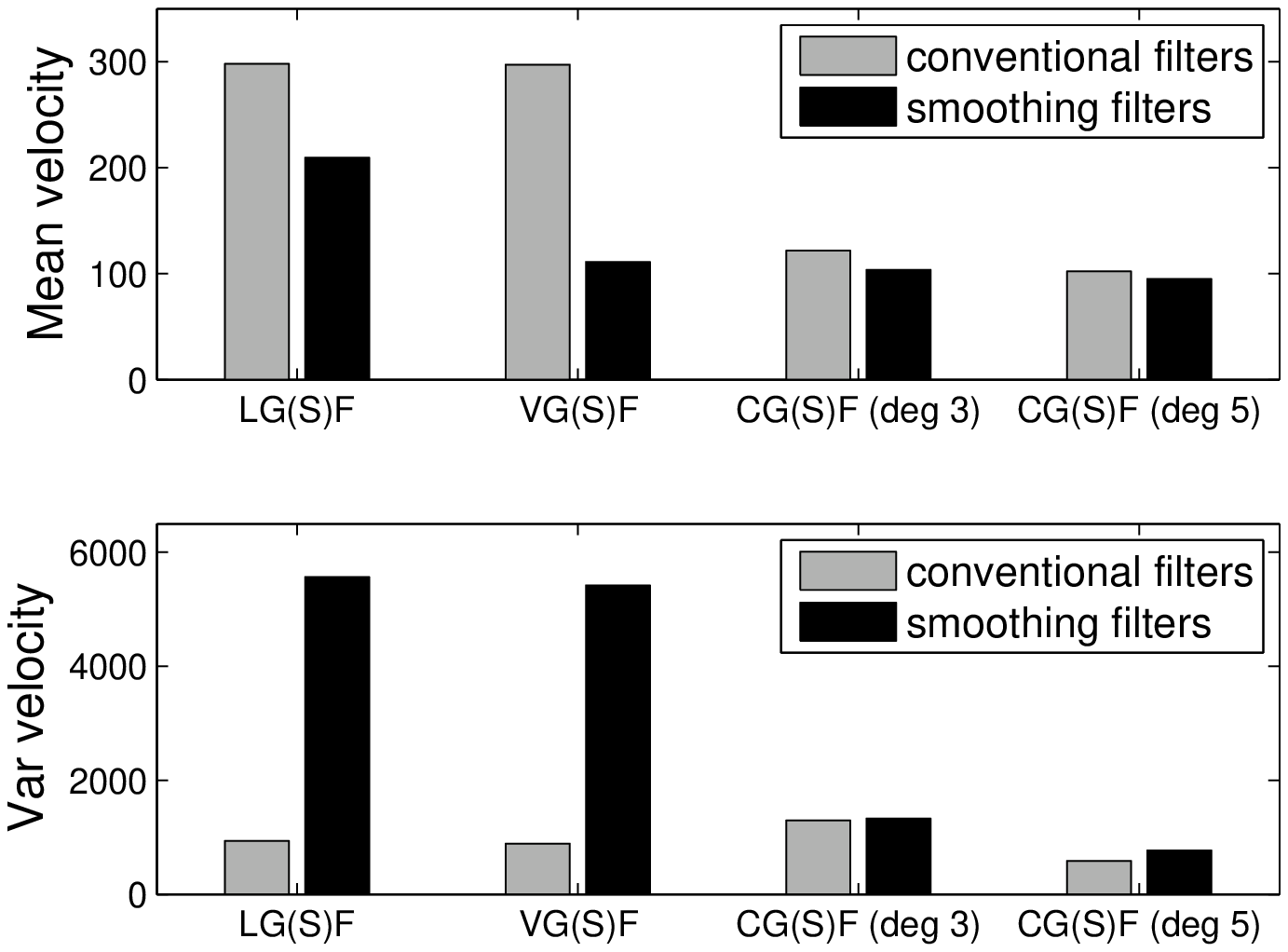} }
\subfigure[turn rate]
{\includegraphics[width=0.43\textwidth]{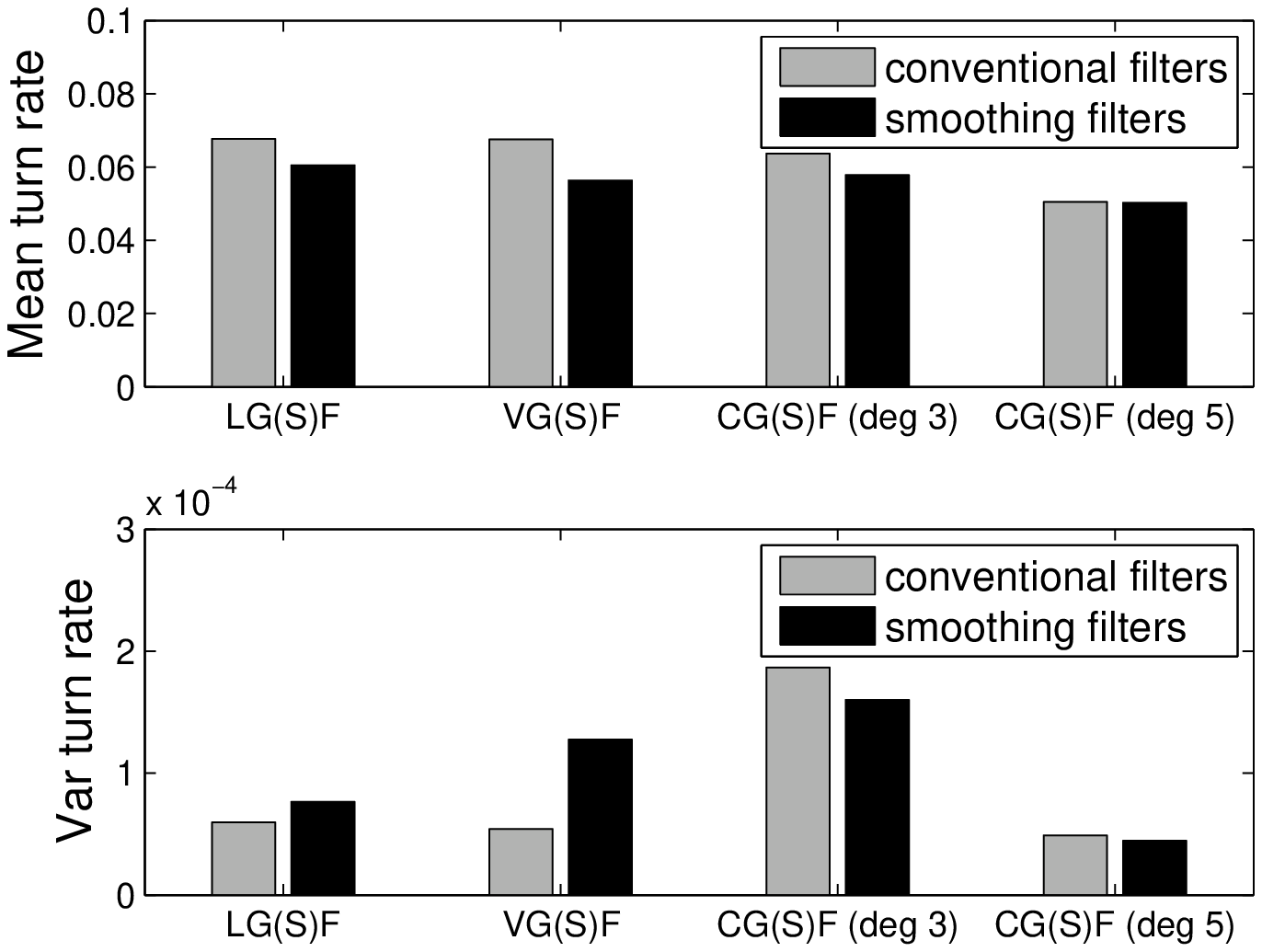} }
\caption{
The mean and variance of time averaged RMSEs between target and filtering estimations.
}
\label{fig:bea4} 
\end{figure}

Fig.~\ref{fig:bea1} displays
one instance of the aircraft trajectory
together with the various filtering estimates.
In this example
the LGF cannot accurately estimate the target 
with the strong nonlinearity
(since the curvature of the position trajectory becomes large near $n=150$)
but the corresponding smoothing filter
never loses the target.
Fig.~\ref{fig:bea2} shows
the average RMSEs,
committed by each filter across $200$ independent simulations,
with respect to position,
velocity
and turn rate.
We see that 
the non-point-based conventional filters 
(LGF and VGF)
become quite in error
around $n=150$, whereas the CGF estimations
keep reasonable accuracy for the entire time period.
The application of smoothing filters provides accuracy improvements in all cases.
To quantify the improvement we turn our attention to the time average.
We depict,
in Fig.~\ref{fig:bea3}, 
the RMSEs averaged over time intervals $50 \leq n \leq 200$.
Although 
the overall accuracy of CGSF is superior to that of LGSF and VGSF,
these two non-point-based smoothing algorithms 
sometimes reach very high accuracy
in the sense of a reduced time average.
Compared with conventional filters,
the enhanced but less uniform accuracy of LGSF and VGSF is illustrated
in terms of the mean and variance
of these temporal RMSEs
in Fig.~\ref{fig:bea4}.

\begin{figure}[h!]
\includegraphics[width=0.43\textwidth]{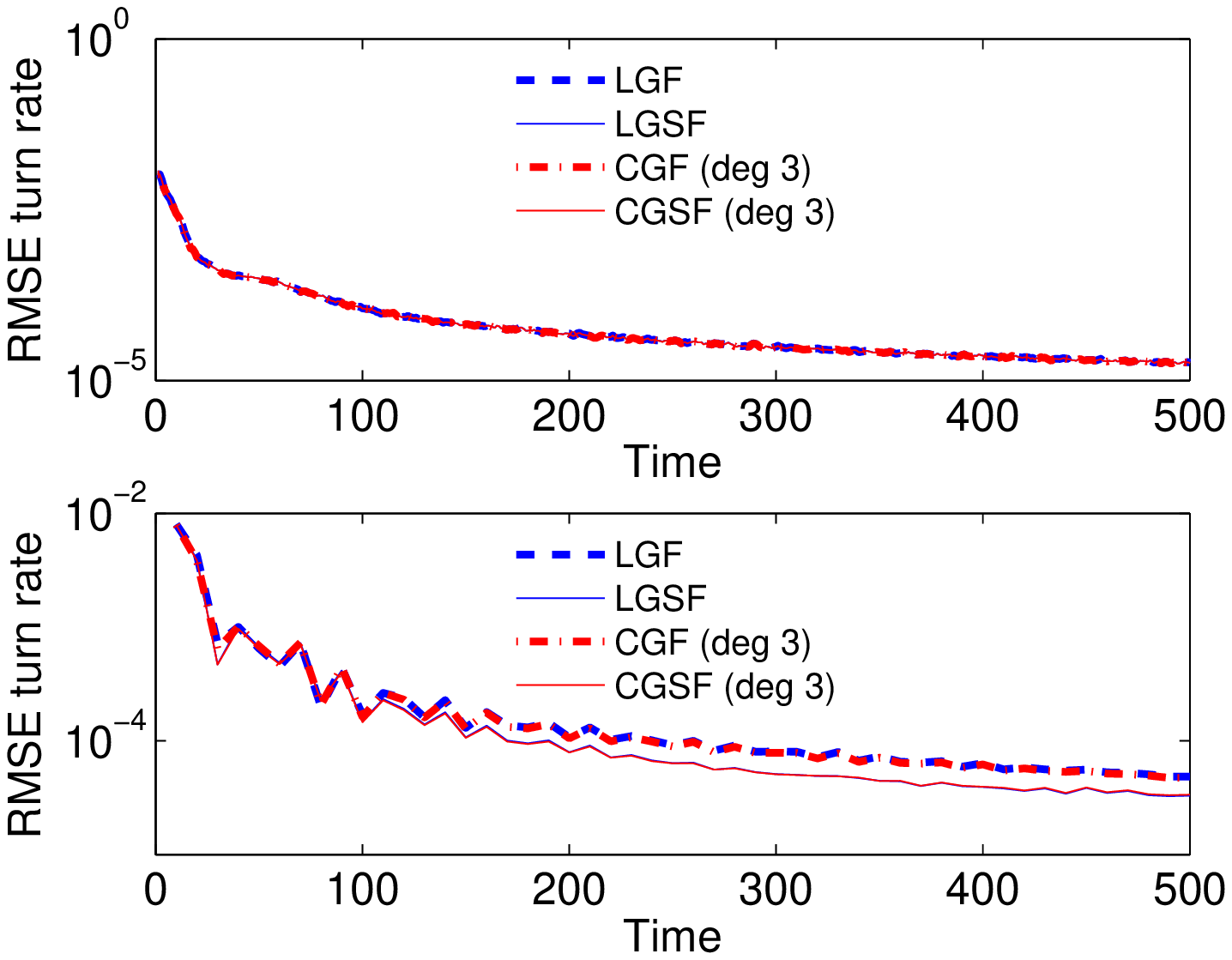}
\caption{
  The RMSEs between the turn rate parameter $\Omega_n$ and its estimates,
  obtained from averaging of $200$ independent simulations.
 The top case uses frequent observations ($\Delta t =\delta t$) 
  and the lower case, sparse observations ($\Delta t =10 \times \delta t$).
}
\label{fig:bea5} 
\end{figure}

Finally we study the system with $q=0$. 
In this case 
$\Omega_n$ is constant
and
the filtering solution can be used for parameter estimation.
In Fig.~\ref{fig:bea5},
we see the smoothing filters outperform conventional filters
particularly with temporally sparse observations.

\section{Conclusion}
\label{sec:conclusion}
This paper formulates a family of sequential Gaussian approximation filters that, 
in contrast to the conventional approaches,
achieve data assimilation via one step backward smoothing
for the solution of the state estimation problem.
The approximate solutions obtained from the proposed smoothing filters
tend to be closer to the observation forward in time
due to the bias of the driving noise conditioned on future observations,
and as a result
can be more accurate than conventional Gaussian filters.
Our numerical simulations, 
performed
on some stochastic systems widely used in the data assimilation community,
show that this is indeed the case
as far as the nonlinearity is involved in either 
the time process equation or the measurement function.
This result is encouraging 
and  leads us to conjecture similar improvements in accuracy
when the smoothing filters are generalised to use Gaussian sum approximations 
in solving the nonlinear filtering problem.

\section*{Acknowledgment}
The authors would like to thank 
King Abdullah University of Science and Technology (KAUST)
Award No. KUK-C1-013-04
for its financial support of this research.

\bibliographystyle{IEEEtran}
\bibliography{smfilter}

%
%
%

\end{document}